\documentclass[11pt]{article}
\usepackage[margin=1in]{geometry}
\usepackage[T1]{fontenc}
\usepackage[utf8]{inputenc}
\usepackage{lmodern}
\usepackage{amsmath,amssymb,amsthm,mathtools}
\usepackage{bm}
\usepackage{physics}
\usepackage{enumitem}
\usepackage{hyperref}
\usepackage{xcolor}
\usepackage{tikz}
\usepackage{microtype}
\usepackage{mathrsfs}
\usepackage{authblk} 
\usepackage[numbers,sort&compress]{natbib} 

\hypersetup{
  colorlinks=true,
  linkcolor=blue,
  citecolor=blue,
  urlcolor=blue
}

\newtheorem{theorem}{Theorem}[section]

\theoremstyle{definition}
\newtheorem{definition}[theorem]{Definition}

\theoremstyle{remark}

\newcommand{\R}{\mathbb{R}}

\newcommand{\diag}{\operatorname{diag}}
\newcommand{\Had}{\circ}

\title{Optimization Workshop Notes for Mathematical Programming with Equilibrium Constraints Algorithms: Penalty Interior-Point, Implicit-Programming, and Piecewise SQP}
\author{Jiguang Yu\thanks{Email: jyu678@bu.edu}}
\affil{College of Engineering, Boston University, Boston, 02215, MA, USA}

\begin{document}
\maketitle

\begin{abstract}

In this workshop, we discuss several algorithms for mathematical programs with equilibrium
constraints (MPECs). The unifying theme is that MPECs are optimization problems whose
feasible set contains a lower-level equilibrium system, often written through complementarity
or variational-inequality conditions. This destroys the smooth manifold or convex structure that standard
nonlinear programming methods rely on. We focus on four algorithmic viewpoints:
\begin{enumerate}[label=(\roman*),itemsep=0.2em]
    \item the classical penalty interior-point algorithm (PIPA);
    \item a monotone-linear complementarity problem (LCP) variant of PIPA that explicitly controls complementarity decay;
    \item an implicit-programming descent method for variational inequality (VI)-constrained MPECs;
    \item piecewise SQP (PSQP), which applies SQP on locally selected smooth pieces.
\end{enumerate}
For each method we explain the model, the search direction subproblem, the globalization
mechanism, and the meaning of the convergence result. Particular emphasis is placed on
what the assumptions are really doing and on the distinction between an attractive
algorithmic idea and a fully valid convergence theorem.

\end{abstract}

\tableofcontents

\section{Introduction}

\subsection{Why MPECs are difficult}

A mathematical program with equilibrium constraints (MPEC) is an optimization problem in
which some constraints are themselves equilibrium conditions. Canonical examples include
complementarity systems, variational inequalities, and KKT systems of lower-level
optimization problems.

A basic abstract form is
\begin{equation}\label{eq:gen-mpec}
\begin{aligned}
\min_{x,y,w,z}\quad & f(x,y,w,z)\\
\text{s.t.}\quad & x \in X,\\
& F(x,y,w,z)=0,\\
& y \ge 0,\quad w \ge 0,\quad y^\top w = 0,
\end{aligned}
\end{equation}
where $f:\mathbb{R}^{n+2m+l}\to \mathbb{R}$ is continuously differentiable on an open set containing $X\times \mathbb{R}_+^{2m}\times \mathbb{R}^l$, \(X\subseteq \R^n\) is typically polyhedral, and $F:\mathbb{R}^{n+2m+l}\to \mathbb{R}^{m+l}$ is twice continuously differentiable on an open set containing $X\times \mathbb{R}^{2m}\times \mathbb{R}^l$.

The difficulty is geometric. The complementarity relation
\[
y_i\ge 0,\qquad w_i\ge 0,\qquad y_iw_i=0
\]
means that each component pair \((y_i,w_i)\) lies on the union of two coordinate axes.
Hence the feasible region is generally nonconvex and assembled from several smooth pieces.
At degenerate indices, where both \(y_i\) and \(w_i\) vanish, the local geometry is not
captured by a single smooth linearization.

\subsection{Three algorithmic philosophies}

The methods studied in these notes correspond to three distinct ways of coping with this
geometry.

\paragraph{Interior-point regularization.}
Maintain \(y>0\) and \(w>0\), replace exact complementarity by a perturbed relation
\(y_iw_i\approx \mu\), and drive \(\mu\downarrow 0\). This yields PIPA.

\paragraph{Reduced-space or implicit programming.}
When the lower-level equilibrium can locally be represented as \(y=y(x)\), optimize the
reduced objective \(x\mapsto f(x,y(x))\).

\paragraph{Piecewise smooth modeling.}
Instead of forcing the feasible region into one global smooth model, select one admissible local
piece and apply SQP there. This yields PSQP.

\subsection{Notation and standing conventions}

Throughout:
\begin{enumerate}[label=(\roman*)]
    \item \(f\) is continuously differentiable unless stated otherwise;
    \item \(X\) or \(\mathcal X\) is closed and convex, often polyhedral;
    \item \(Y\coloneqq\diag(y)\) and \(W\coloneqq\diag(w)\);
    \item \(\mu \coloneqq \frac{y^\top w}{m}\) denotes the average complementarity product.
\end{enumerate}

\section{Stationarity for complementarity-constrained MPECs}

Consider \eqref{eq:gen-mpec}. Let
\[
\mathcal F \coloneqq
\{(x,y,w,z): x\in X,\ F(x,y,w,z)=0,\ y\ge 0,\ w\ge 0,\ y^\top w=0\}.
\]
be the feasible region.

\subsection{Index sets and strict complementarity}

At a feasible point \(u^*=(x^*,y^*,w^*,z^*)\), define
\[
\alpha\coloneqq\{i:y_i^*>0=w_i^*\},\qquad
\beta \coloneqq\{i:y_i^*=0=w_i^*\},\qquad
\gamma\coloneqq\{i:y_i^*=0<w_i^*\}.
\]
The set \(\beta\) is the degenerate set.

\begin{definition}[Strict complementarity]
A feasible point is strictly complementary if
\[
y^*+w^*>0
\quad\text{componentwise.}
\]
Equivalently, \(\beta=\varnothing\).
\end{definition}

Strict complementarity is not merely a cosmetic simplification. It removes the
branching ambiguity of the local feasible geometry.

\subsection{First-order stationarity}

A feasible point \(u^*\in\mathcal F\) is stationary if
\[
\nabla f(u^*)^\top d \ge 0
\qquad \forall d\in T(u^*;\mathcal F),
\]
where \(T(u^*;\mathcal F)\) is the tangent cone.

Under strict complementarity, one can express stationarity through a KKT-like system.
Let
\[
M^*=
\begin{bmatrix}
\nabla_yF(u^*) & \nabla_wF(u^*) & \nabla_zF(u^*)\\
W^* & Y^* & 0
\end{bmatrix},
\qquad
Y^*\coloneqq\diag(y^*),\quad W^*\coloneqq\diag(w^*).
\]

\begin{theorem}[KKT form under strict complementarity]
Assume \(u^*\) is feasible, strictly complementary, and \(M^*\) is nonsingular.
Then \(u^*\) is stationary if and only if there exist multipliers
\[
\zeta,\ \pi,\ \xi
\]
such that
\begin{align*}
\nabla_x f(u^*) - \nabla_xF(u^*)^\top \pi + G^\top \zeta &= 0,\\
\nabla_y f(u^*) - \nabla_yF(u^*)^\top \pi - W^*\xi &= 0,\\
\nabla_w f(u^*) - \nabla_wF(u^*)^\top \pi - Y^*\xi &= 0,\\
\nabla_z f(u^*) - \nabla_zF(u^*)^\top \pi &= 0,\\
\zeta\ge 0,\qquad \zeta^\top(a-Gx^*) &= 0.
\end{align*}
\end{theorem}

The additional multiplier \(\xi\) records the effect of complementarity.

\section{The classical penalty interior-point algorithm (PIPA)}

\subsection{Merit and infeasibility measures}

For \((x,y,w,z)\) with \(y,w\ge 0\), define
\[
\phi(x,y,w,z)\coloneqq\lVert F(x,y,w,z)\rVert^2 + y^\top w
\]
and, for \(\alpha>0\),
\[
P_\alpha(x,y,w,z)\coloneqq f(x,y,w,z)+\alpha\,\phi(x,y,w,z).
\]

Thus \(\phi=0\) is equivalent to lower-level feasibility plus complementarity, while
\(P_\alpha\) balances objective reduction against movement toward feasibility.

\subsection{Interior-point geometry}

PIPA keeps all iterates strictly positive in the complementarity variables:
\[
y^v>0,\qquad w^v>0.
\]
It also imposes the centrality condition
\[
y_i^v w_i^v \ge p\,\mu_v,
\qquad
\mu_v\coloneqq\frac{(y^v)^\top w^v}{m},
\]
for some fixed \(p\in(0,1)\).

The point of centrality is to prevent one pair \((y_i,w_i)\) from collapsing much faster
than the others. This is analogous to staying in a neighborhood of the central path.

\subsection{Direction-finding subproblem}

At the current iterate \(u^v=(x^v,y^v,w^v,z^v)\), define
\[
Y^v\coloneqq\diag(y^v),\qquad W^v\coloneqq\diag(w^v),\qquad F^v\coloneqq F(u^v).
\]
A search direction \(d^v=(dx^v,dy^v,dw^v,dz^v)\) is computed from the quadratic subproblem
\begin{equation}\label{eq:pipa-qp}
\begin{aligned}
\min_{dx,dy,dw,dz}\quad &
\nabla f(u^v)^\top
\begin{bmatrix}
dx\\dy\\dw\\dz
\end{bmatrix}
+\frac12 dx^\top Q^v dx\\
\text{s.t.}\quad &
x^v+dx\in X,\\
& \lVert dx\rVert^2 \le c\bigl(\lVert F^v\rVert+(y^v)^\top w^v\bigr),\\
& \nabla_xF^v dx + \nabla_yF^v dy + \nabla_wF^v dw + \nabla_zF^v dz = -F^v,\\
& W^vdy + Y^v dw = -Y^vw^v + \sigma_v\mu_v e.
\end{aligned}
\end{equation}

This has three coupled components:
\begin{enumerate}[label=(\roman*),itemsep=0.2em]
    \item an SQP model for the upper-level objective in \(x\);
    \item a Newton step for the lower-level equation \(F=0\);
    \item a primal-dual centering equation for complementarity.
\end{enumerate}

The unusual trust-region-like restriction
\begin{equation}\label{eq:stepp_size}
\lVert dx\rVert^2 \le c\bigl(\lVert F^v\rVert+(y^v)^\top w^v\bigr)
\end{equation}
links progress in \(x\) to lower-level infeasibility. This is analytically convenient but
also the source of an important pathology: as feasibility improves, the allowed \(x\)-step
may become artificially tiny.

\subsection{Well-posedness mechanism}

The lower-level linearized block
\[
\begin{bmatrix}
\nabla_yF & \nabla_wF & \nabla_zF\\
W & Y & 0
\end{bmatrix}
\]
must be nonsingular whenever \(y>0\) and \(w>0\). In the standard theory this is implied by
a mixed \(P_0\)-type assumption on the Jacobian of the lower system. The operational meaning
is simple: once \(dx\) is chosen, \((dy,dw,dz)\) are uniquely determined.

\subsection{Line search and penalty update}

The trial arc is
\[
u^v(\tau)=u^v+\tau d^v,\qquad \tau\in[0,1].
\]
The line search enforces:
\begin{enumerate}[label=(\roman*),itemsep=0.2em]
    \item positivity of \(y^v(\tau)\) and \(w^v(\tau)\);
    \item preservation of centrality;
    \item sufficient decrease in \(\phi\) and in the penalty merit \(P_{\alpha_v}\).
\end{enumerate}

Under boundedness of a suitable level set and the mixed-\(P_0\) assumption, the classical
analysis shows that iterates remain bounded, stay strictly interior, preserve centrality,
and possess accumulation points.

Under additional assumptions at a cluster point --- notably strict complementarity and
nonsingularity of the Jacobian block \(M^*\) --- the intended conclusion is stationarity.

However, the strongest global convergence claim is not correct in full generality.
Leyffer constructed a counterexample showing that PIPA can converge to a nonstationary
point. The mechanism is precisely that the step restriction \eqref{eq:stepp_size} on \(dx\) can suppress motion
in the upper-level variable too aggressively as complementarity and lower-level infeasibility
become small. One should therefore retain the algorithmic insight of PIPA while treating
the original theorem with caution.

\section{A monotone-LCP variant of PIPA}

\subsection{Problem class}

A particularly transparent special case is the monotone linear-complementarity lower level:
\begin{equation}\label{eq:lcp-mpec}
\begin{aligned}
\min_{x,y,w}\quad & f(x,y)\\
\text{s.t.}\quad & x\in X\coloneqq\{x:Gx\le a\},\\
& w = q + Nx + My,\\
& y\ge 0,\quad w\ge 0,\quad y^\top w = 0,
\end{aligned}
\end{equation}
where \(M\succeq 0\).

Define the affine residual
\[
r(x,y,w)\coloneqq q+Nx+My-w
\]
and infeasibility measure
\[
\phi(x,y,w)\coloneqq y^\top w + \lVert r(x,y,w)\rVert.
\]

\subsection{Direction subproblem}

At iteration \(v\), with \(r^v=r(x^v,y^v,w^v)\), the direction \((dx^v,dy^v,dw^v)\) solves
\begin{equation}\label{eq:lcp-qp}
\begin{aligned}
\min_{dx,dy,dw}\quad &
(\nabla_x f^v)^\top dx + (\nabla_y f^v)^\top dy + \frac12 dx^\top Q_v dx\\
\text{s.t.}\quad &
x^v+dx\in X,\\
& W_vdy + Y_vdw = -Y_vw^v + \sigma_v\mu^v e,\\
& Ndx + Mdy - dw = -r^v,\\
& \lVert dx\rVert ^2 \le c\,\phi^v.
\end{aligned}
\end{equation}

Because the residual equation is affine, the trial arc
\[
x^v(r)=x^v+r\,dx^v,\qquad y^v(r)=y^v+r\,dy^v,\qquad w^v(r)=w^v+r\,dw^v
\]
satisfies the exact identity
\[
r(x^v(r),y^v(r),w^v(r))=(1-r)\,r^v.
\]
Thus the affine feasibility residual decays linearly with \(1-r\).

\subsection{Limited complementarity decrease}

The main modification relative to classical PIPA is the additional condition
\[
(w^v(r))^\top y^v(r)\ge (1-r)\,(w^v)^\top y^v.
\]
By this requirement, a higher priority to satisfying the linear constraint $q+Nx+My-w=0$ is given than to the complementary condition $y^\top w=0$; especially, complementarity is expected not to be satisfied before feasibility.
This is an elegant repair of the imbalance that can arise in the original PIPA logic.

\subsection{Why the step-size analysis is explicit}

Along the arc,
\[
w(r)\Had y(r)=w\Had y + r(w\Had dy+y\Had dw)+r^2(dw\Had dy),
\]
and the perturbed complementarity equation implies
\[
w\Had dy + y\Had dw = -w\Had y + \sigma\mu e.
\]
Hence
\[
w(r)\Had y(r)
=
(1-r)\,w\Had y + r\,\sigma\mu e + r^2(dw\Had dy).
\]
Summing components gives
\[
w(r)^\top y(r)
=
(1-r)\,w^\top y + \sigma r\,w^\top y + r^2(dw)^\top dy.
\]

Thus both the centrality condition and the limited-complementarity condition reduce to
scalar quadratic inequalities in \(r\). This makes the line search unusually explicit.

\subsection{Convergence message}

The global theorem for this variant has the familiar penalty-method dichotomy:
\begin{enumerate}[label=(\roman*),itemsep=0.2em]
    \item if the penalty sequence blows up, one extracts stationarity information along the
    penalty-update subsequence;
    \item if the penalty stays bounded and accepted step lengths stay away from zero, then
    regular accumulation points are stationary.
\end{enumerate}

Compared with classical PIPA, the conceptual novelty is not the QP itself, but the
line-search rule that keeps complementarity reduction synchronized with residual reduction.

\section{Implicit-programming descent for VI-constrained MPECs}

\subsection{The reduced-space viewpoint}

Now consider MPECs of the form
\begin{equation}\label{eq:vi-mpec}
\begin{aligned}
\min_{x,y}\quad & f(x,y)\\
\text{s.t.}\quad & (x,y)\in \mathcal F \coloneqq (\mathcal X\times \R^m)\cap \operatorname{Gr}\mathcal S,
\end{aligned}
\end{equation}
where \(\mathcal S(x)\) is the solution set of a parametric variational inequality.

The central idea is to assume that locally one can select the lower-level solution by an
implicit map \(y(\cdot)\). Then the MPEC is treated as the reduced problem
\[
\min_{x\in \mathcal X}\ f(x,y(x)).
\]

\subsection{Local implicit solution map}

At a feasible point \((x,y)\), suppose there exists a local Lipschitz mapping
\[
\hat y:V\cap \mathcal X\to U
\]
such that \(\hat y(x)=y\) and \(\hat y(\xi)\in \mathcal S(\xi)\) for nearby \(\xi\).
Assume further that \(\hat y\) is directionally differentiable at \(x\).
Then the reduced directional derivative is
\[
\varphi'(x;dx)
=
\nabla_x f(x,y)^\top dx
+
\nabla_y f(x,y)^\top \hat y'(x;dx).
\]
We thus have the first-order stationarity for the reduced problem.
\begin{definition}[Stationarity via the implicit map]
A feasible point \((x,y)\) is stationary if
\[
\nabla_x f(x,y)^\top dx
+
\nabla_y f(x,y)^\top \hat y'(x;dx)
\ge 0
\qquad \forall dx\in T(x;\mathcal X).
\]
\end{definition}

\subsection{Direction-finding model}

At iterate \((x^\nu,y^\nu)\), with local implicit map \(y^\nu(\cdot)\), define the model
subproblem
\begin{equation}\label{eq:implicit-subproblem}
\begin{aligned}
\min_{dx}\quad &
(\nabla_x f^\nu)^\top dx
+
(\nabla_y f^\nu)^\top (y^\nu)'(x^\nu;dx)
+
\frac12 dx^\top Q_\nu dx\\
\text{s.t.}\quad & x^\nu+dx\in \mathcal X.
\end{aligned}
\end{equation}

This is a regularized first-order model of the reduced objective. Even if \(Q_\nu\succ 0\),
the problem need not be convex, because \(dx\mapsto (y^\nu)'(x^\nu;dx)\) may be nonlinear.

\subsection{Line search along the equilibrium path}

The update is not along a straight line in \((x,y)\)-space. Instead one uses
\[
r\mapsto \bigl(x^\nu+r\,dx^\nu,\ y^\nu(x^\nu+r\,dx^\nu)\bigr).
\]
Armijo backtracking is performed on this equilibrium path.

This is the conceptual hallmark of the method: lower-level feasibility is preserved by
construction, because the path remains on the graph of the solution map.

\subsection{Convergence logic}

The convergence theorem typically requires:
\begin{enumerate}[label=(\roman*),itemsep=0.2em]
    \item polyhedrality of \(\mathcal X\);
    \item closedness of the solution map;
    \item uniform positive definiteness and boundedness of \(Q_\nu\);
    \item boundedness of a lower level set of \(f\);
    \item strong Fr\'echet differentiability of the local implicit map at the cluster point.
\end{enumerate}

Under these assumptions, any cluster point is stationary.

\subsection{Main practical limitation}

The reduced-space formulation is elegant, but computation can be difficult:
\begin{enumerate}[label=(\roman*),itemsep=0.2em]
    \item the derivative \(y'(x;dx)\) may itself require solving an auxiliary VI or
    complementarity problem;
    \item each line-search trial may require a fresh lower-level solve.
\end{enumerate}
So the theory is conceptually clean, but the method is expensive unless the implicit map
has a tractable derivative.

\section{Piecewise SQP (PSQP)}

SQP works extremely well for smooth nonlinear programs, but direct SQP on MPECs is
fragile near degenerate complementarity points. The reason is that one linearization cannot
faithfully represent a feasible region that is locally a union of several smooth pieces.

PSQP resolves this by applying SQP on one selected local branch at a time.

\subsection{MPEC in KKT form}

Consider an MPEC written as
\begin{equation}\label{eq:kkt-mpec}
\begin{aligned}
\min_{x,y,\lambda}\quad & f(x,y)\\
\text{s.t.}\quad & Gx+Hy+a\le 0,\\
& \lambda\ge 0,\\
& g(x,y)\le 0,\\
& L(x,y,\lambda)=0,\\
& \lambda^\top g(x,y)=0,
\end{aligned}
\end{equation}
where
\[
L(x,y,\lambda)=F(x,y)+\nabla_y g(x,y)^\top\lambda.
\]

At a feasible point define
\[
I_0(z,\lambda)\coloneqq\{i:g_i(z)=0=\lambda_i\},
\qquad
I_+(z,\lambda)\coloneqq\{i:g_i(z)=0<\lambda_i\}.
\]
The degenerate set \(I_0\) is exactly where multiple local pieces arise.

\subsection{Piecewise NLPs}

For a partition \(J_1\cup J_2=\{1,\dots,\ell\}\) consistent with the current complementarity
regime, define the local NLP piece
\begin{equation}\label{eq:piece-nlp}
\begin{aligned}
\min_{x,y,\lambda}\quad & f(x,y)\\
\text{s.t.}\quad & Gx+Hy+a\le 0,\\
& L(x,y,\lambda)=0,\\
& \lambda_i=0,\quad g_i(x,y)\le 0,\qquad i\in J_1,\\
& \lambda_i\ge 0,\quad g_i(x,y)=0,\qquad i\in J_2.
\end{aligned}
\end{equation}

Thus each degenerate index is assigned to one branch or the other.

\subsection{PSQP subproblem}

At the current iterate \(w^v=(x^v,y^v,\lambda^v)\), with multiplier estimate \(\nu^v\), and
for one admissible partition \((J_1,J_2)\), solve
\begin{equation}\label{eq:psqp}
\begin{aligned}
\min_{dw}\quad &
\nabla f(z^v)^\top dz
+\frac12 dw^\top \nabla_{ww}^2\mathcal L_{\mathrm{MPEC}}(w^v,\nu^v)\,dw\\
\text{s.t.}\quad &
G(x^v+dx)+H(y^v+dy)+a\le 0,\\
& L(z^v,\lambda^v)+\nabla L(z^v,\lambda^v)\,dw=0,\\
& (\lambda^v+d_\lambda)_i=0,\quad g_i(z^v)+\nabla g_i(z^v)^\top dz\le 0,
\quad i\in J_1,\\
& (\lambda^v+d_\lambda)_i\ge 0,\quad g_i(z^v)+\nabla g_i(z^v)^\top dz=0,
\quad i\in J_2.
\end{aligned}
\end{equation}

This is just an SQP step --- but on one selected smooth piece.

\subsection{Why local convergence works}

The local convergence theorem has a finite-piece argument:
\begin{enumerate}[label=(\roman*),itemsep=0.2em]
    \item for each admissible local piece, ordinary SQP theory applies;
    \item only finitely many pieces can arise near a fixed solution;
    \item hence the local contraction estimates can be uniformized across pieces.
\end{enumerate}

The result is local Q-superlinear convergence, and often Q-quadratic convergence under
stronger smoothness assumptions.

\subsection{What PSQP teaches conceptually}

PSQP shows that one need not regularize complementarity away in order to justify
Newton/SQP-type local behavior. One may instead respect the piecewise geometry directly.

\section{Comparison of the four viewpoints}

\subsection{What each method linearizes}

\begin{center}
\begin{tabular}{p{3.2cm}p{9.3cm}}
\hline
Method & Local model idea\\
\hline
PIPA & Linearize lower equations and perturbed complementarity while keeping \(y,w>0\).\\
LCP-PIPA variant & Same, but synchronize complementarity decay with affine residual decay.\\
Implicit programming & Eliminate the lower-level variable locally through \(y=y(x)\).\\
PSQP & Pick one local complementarity branch and apply ordinary SQP there.\\
\hline
\end{tabular}
\end{center}

\subsection{Strengths and weaknesses}

\paragraph{PIPA.}
Strongly structured, interior-point flavored, and conceptually appealing; but the original
global convergence theory is too optimistic.

\paragraph{Monotone-LCP variant.}
Cleaner balance between residual reduction and complementarity reduction; best viewed as a
targeted repair of classical PIPA for an important special case.

\paragraph{Implicit programming.}
Theoretically elegant and faithful to the lower-level equilibrium map; computationally hard
because the implicit derivative and line-search evaluations may be expensive.

\paragraph{PSQP.}
Excellent local interpretation and fast local rates under strong assumptions; but by itself it
is a local method and does not solve the globalization problem.

\newpage
\section{Practice Questions}
\subsection{True or False}

\begin{enumerate}[label=(\roman*)]
    \item In the primal-dual formulation of the stationary conditions for an MPEC, the piecewise linear map $LH^*$ being a global Lipschitzian homeomorphism is equivalent to the pair of matrices $(W_0, W_1)$ possessing the W property.

\item If the partitioned matrix $Q = [A \quad B \quad C]$ has the mixed P property, then for any subset $\alpha$, the submatrix consisting of the columns $A_\alpha$, the complement columns $B_{\overline{\alpha}}$, and all columns of $C$ must be positive definite.

\item The convergence analysis for the general Penalty Interior Point Algorithm (PIPA) relies on the strict complementarity assumption ($y^* + w^* > 0$) at the limit point, a requirement that cannot currently be removed for MPECs, despite being inessential for standard complementarity problems.

\item In the alternative PIPA designed for LCP-constrained mathematical programs, the ``limited complementarity decrease'' condition is enforced to ensure that the complementarity condition ($y^Tw = 0$) is satisfied strictly before the linear feasibility constraints.

\item Under the Global B-differentiability and Implicit Function (GBIF) assumption, the direction-finding subproblem in the implicit programming descent algorithm is guaranteed to be a strictly convex quadratic program because the implicit function is Lipschitz continuous.

\item In the implicit programming descent algorithm, if the globally optimum objective value of the direction-finding subproblem is strictly negative, the current iterate cannot be a stationary point of the MPEC.

\item A direct application of the standard Sequential Quadratic Programming (SQP) method to an MPEC formulated with KKT constraints often fails at degenerate vectors because the Strict Mangasarian-Fromovitz Constraint Qualification (SMFCQ) generally requires strict complementarity to hold.

\item The Piecewise Sequential Quadratic Programming (PSQP) algorithm overcomes the limitations of standard SQP by enforcing a strict complementarity assumption at the limit point to ensure superlinear convergence.

\item In the convergence proof for the implicit programming descent algorithm, the tangent cone property for polyhedra guarantees that for any direction $d$ in the tangent cone at a point $\overline{x}$, taking a sufficiently small step along $d$ from any feasible point $y$ sufficiently close to $\overline{x}$ will remain within the polyhedron.

\item The PSQP method for MPECs and the Kojima-Shindo method for piecewise smooth equations are identical in their iterative subproblems; both rely on solving a system of linear equations obtained from a single Newton iteration on an active piece.

\end{enumerate}

\subsection{Computational questions}

\subsubsection*{Problem 1: The $LH^*$ Map and Strict Complementarity (30 Points)}

Consider a simplified MPEC where the inner problem is formulated as a parametric system of equations and complementarity conditions. Let $x \in \mathbb{R}$ be the upper-level variable and $(y, w) \in \mathbb{R}^2$ be the lower-level variables. The equality constraint $F(x,y,w) = 0$ is given by:
\[
F(x,y,w) = x + y - w = 0
\]
subject to the complementarity condition $(y, w) \ge 0, \; y w = 0$.

Let the current point be evaluated at $u^* = (x^*, y^*, w^*) = (1, 0, 0)$.

\begin{enumerate}
    \item[(a)] Evaluate the index sets $\alpha, \beta,$ and $\gamma$ for the lower-level variables at $u^*$ as defined in Section 6.1.1. Is this point degenerate?
    \item[(b)] Construct the piecewise linear mapping $LH^*(dy, dw)$ at $u^*$ for the relevant variables.
    \item[(c)] Using Lemma 6.1.1, determine whether the map $LH^*$ is a global Lipschitzian homeomorphism. Show your work by explicitly constructing the matrix pair $(W_0, W_1)$ and testing the $\mathcal{W}$ property.
\end{enumerate}

\subsubsection*{Problem 2: PIPA Quadratic Subproblem Setup (35 Points)}

Consider the following LCP-constrained Mathematical Program (Section 6.2):
\begin{align*}
    \text{minimize} \quad & f(x,y) = \frac{1}{2}x^2 + xy + y^2 \\
    \text{subject to} \quad & x \ge 0 \\
    & w = -2 + x + y \\
    & (y, w) \ge 0, \quad y w = 0
\end{align*}

You are applying the variant of PIPA for LCP-constrained MPs. Assume the current strictly interior iterate is $x^0 = 1, y^0 = 2, w^0 = 1$. 

Let the algorithmic parameters be:
\begin{itemize}
    \item Centering parameter: $\sigma_0 = 0.5$
    \item Step-bound parameter: $c = 10$
    \item Hessian approximation: $Q_0 = [1]$ (a $1 \times 1$ matrix for the $dx$ term)
\end{itemize}

\begin{enumerate}
    \item[(a)] Calculate the current penalty constraint residual $\phi_0 = y^0 w^0 + \|r^0\|$, where $r^0$ is the violation of the linear equality $w = q + Nx + My$. 
    \item[(b)] Explicitly write down the full Quadratic Program ($QP_0$) in the variables $(dx, dy, dw)$ that must be solved to find the search direction. Plug in all numerical values for gradients, matrices, and state variables. \textit{(You do not need to solve the QP).}
\end{enumerate}

\subsubsection*{Problem 3: Implicit Programming and Piecewise Subproblems (35 Points)}

Consider an MPEC defined over a closed convex set $X = [-1, 1]$:
\begin{align*}
    \text{minimize} \quad & f(x,y) = x^2 - y \\
    \text{subject to} \quad & x \in [-1, 1] \\
    & y \in \text{SOL}(F(x, \cdot), \mathbb{R}_+)
\end{align*}
where $F(x,y) = y - x$. The lower-level problem is equivalent to the LCP: $y \ge 0, \; y - x \ge 0, \; y(y-x) = 0$.

Assume the current iterate is $x^\nu = 0$, and the corresponding optimal lower-level response is $\overline{y}^\nu = 0$. Let the Hessian approximation be $Q_\nu = [2]$.

\begin{enumerate}
    \item[(a)] The solution mapping $\overline{y}(x)$ is B-differentiable but not Fréchet differentiable at $x^\nu = 0$. Calculate the directional derivative $d\overline{y} = (\overline{y}^\nu)'(0; dx)$ as a piecewise function of the direction $dx \in \mathbb{R}$.
    \item[(b)] Formulate the implicit programming descent subproblem (Equation 6.3.2) in the single variable $dx$. 
    \item[(c)] Solve the subproblem from part (b) by hand. Find the globally optimal search direction $dx^*$. Does this direction yield a strictly negative objective value in the subproblem, and what does Proposition 6.3.2 imply about the stationarity of $x^\nu = 0$?
\end{enumerate}

\subsection{Proof questions}

\subsubsection*{Problem 1: Matrix Properties in Complementarity}
In the context of the Penalty Interior Point Algorithm (PIPA) for solving an MPEC, the properties of the partitioned Jacobian matrix are critical for ensuring the non-singularity of the search-direction system. Let $Q = \begin{bmatrix} A & B & C \end{bmatrix}$, where $A, B \in \mathbb{R}^{(m+l) \times m}$ and $C \in \mathbb{R}^{(m+l) \times l}$. 

\begin{enumerate}[label=(\alph*)]
    \item Define the \textit{mixed P property} for the partitioned matrix $Q$.
    \item Assume $Q$ possesses the mixed P property. Prove rigorously that for every index subset $\alpha \subseteq \{1,\dots,m\}$ with complement $\overline{\alpha}$, the submatrix constructed as $\begin{bmatrix} A_\alpha & B_{\overline{\alpha}} & C \end{bmatrix}$ has linearly independent columns.
\end{enumerate}

\subsubsection*{Problem 2: PIPA Merit Function and Descent}
Consider the nonnegative constraint violation function $\phi(x,y,w,z) \equiv F(x,y,w,z)^T F(x,y,w,z) + w^T y$ and the corresponding penalized objective function $P_\alpha(x,y,w,z) \equiv f(x,y,w,z) + \alpha\phi(x,y,w,z)$.

Let $(dx^\nu, dy^\nu, dw^\nu, dz^\nu)$ be the unique optimal solution to the QP subproblem at iteration $\nu$, which enforces the perturbed complementarity condition:
\[
w_j dy_j + y_j dw_j = -w_j y_j + \sigma_\nu \frac{y^T w}{m}, \quad \forall j = 1, \dots, m.
\]

\begin{enumerate}[label=(\alph*)]
    \item Derive the exact expression for the directional derivative components of $\phi$ along the Newton direction. Specifically, prove that:
    \[
    (d\phi_x^\nu)^T dx^\nu + (d\phi_y^\nu)^T dy^\nu + (d\phi_w^\nu)^T dw^\nu + (d\phi_z^\nu)^T dz^\nu = -2\|F^\nu\|^2 - (1 - \sigma_\nu)(y^\nu)^T w^\nu.
    \]
    Show all steps of the chain rule and substitute the appropriate QP constraints.
\end{enumerate}

\subsubsection*{Problem 3: Implicit Programming and Stationarity}
Consider the implicit programming formulation of an MPEC where the lower-level problem is parameterized by $x \in X$, and let $\overline{y}(x)$ be the locally Lipschitz, B-differentiable implicit solution function. To find a descent direction, we solve the following constrained minimization problem in $dx \in \mathbb{R}^n$:
\begin{align*}
    \text{minimize} \quad & (df_x^\nu)^T dx + (df_y^\nu)^T (\overline{y}^\nu)'(x^\nu; dx) + \frac{1}{2} dx^T Q_\nu dx \\
    \text{subject to} \quad & x^\nu + dx \in X
\end{align*}
where $X$ is a polyhedron and $Q_\nu$ is a symmetric positive definite matrix.

\begin{enumerate}[label=(\alph*)]
    \item Prove that the optimum objective value of this subproblem is zero if and only if $dx = 0$ is a globally optimal solution.
    \item Prove that if the optimum objective value of the subproblem is zero, then the current iterate $(x^\nu, y^\nu)$ is a stationary point of the implicit MPEC formulation. State any necessary theorems regarding the tangent cone $\mathcal{T}(x^\nu; X)$ required to complete your proof.
\end{enumerate}

\subsubsection*{Problem 4: Polyhedral Geometry in Piecewise SQP}
In the convergence analysis of the Piecewise Sequential Quadratic Programming (PSQP) algorithm and descent methods, the geometric properties of the upper-level constraints dictate step-size validity.

\begin{enumerate}[label=(\alph*)]
    \item Let $X = \{x \in \mathbb{R}^n : Gx \le a\}$ be a polyhedron and let $\overline{x} \in X$. Let $\mathcal{T}(\overline{x}; X)$ denote the tangent cone to $X$ at $\overline{x}$. Prove that if $d \in \mathcal{T}(\overline{x}; X)$, there exist positive scalars $\overline{\epsilon}$ and $\delta$ such that:
    \[
    y + \epsilon d \in X
    \]
    for every $\epsilon \in [0, \overline{\epsilon}]$ and every $y \in X$ such that $\|y - \overline{x}\| \le \delta$.
    \item Briefly explain (1-2 paragraphs) how this geometric property is rigorously applied to execute the limit transition in the convergence proof of the descent algorithm, specifically when analyzing an infinite sequence of iterates approaching a limit point $x^*$.
\end{enumerate}

\section{Solutions for Practice Questions}
\subsection{Solutions for True or False}

\begin{enumerate}[label=(\roman*)]
    \item Answer: True. Lemma 6.1.1 explicitly states that the map $LH^*$ being a global Lipschitzian homeomorphism is equivalent to the pair $(W_0, W_1)$ having the W property.
    \item Answer: False. According to Proposition 6.1.5, the mixed P property implies that this specific submatrix has linearly independent columns, not necessarily that it is positive definite.
    \item Answer: True. The text explicitly notes that the assumption of strict complementarity at the limit point is a present deficiency in the PIPA convergence results for MPECs, even though it is inessential when applying interior point algorithms to standard complementarity problems.
    \item Answer: False. The rationale for the limited complementarity decrease is the exact opposite; it gives higher priority to satisfying the linear constraints, ensuring that complementarity is not satisfied before feasibility.
    \item Answer: False. The objective function of the subproblem is not necessarily strictly convex—or even convex—due to the possible lack of convexity of the directional derivative function $(\overline{y}^\nu)'(x^\nu; \cdot)$. It only becomes a strictly convex quadratic program if the implicit function is Fréchet differentiable and the constraint set $X$ is polyhedral.
    \item Answer: True. Proposition 6.3.2 establishes that a globally optimum objective value of exactly zero is equivalent to the current point being a stationary point of the MPEC. A negative value implies a descent direction exists.
    \item Answer: True. The text notes that for the tangent cone and linearized cone to coincide (implied by SMFCQ), the multiplier must generally be strictly complementary, causing the standard SQP method to frequently fail at degenerate vectors.
    \item Answer: False. The primary advantage of the PSQP algorithm is that it modifies the SQP approach by exploiting the piecewise structure of the feasible set, meaning its convergence does not depend on the strict complementarity assumption.
    \item Answer: True. Lemma 6.3.3 formally states that if $d \in \mathcal{T}(\overline{x}; X)$, there exist positive scalars $\overline{\epsilon}$ and $\delta$ such that $y + \epsilon d \in X$ for every $\epsilon \in [0, \overline{\epsilon}]$ and every $y \in X$ within $\delta$ distance of $\overline{x}$.
    \item Answer: False. While they share the core idea of solving a linearized subproblem on an active piece, they differ primarily because the PSQP subproblem is a quadratic program, whereas the Kojima-Shindo subproblem is a system of linear equations.
\end{enumerate}

\subsection{Solutions to computational questions}

\subsubsection*{Problem 1: The $LH^*$ Map and Strict Complementarity} 
    
(a) Index Sets and Degeneracy

At the given point $u^* = (x^*, y^*, w^*) = (1, 0, 0)$, we evaluate the index sets based on the lower-level variables $(y, w)$:
\begin{align*}
    \alpha &\equiv \{i : y^*_i > 0 = w^*_i\} = \emptyset, \\[4pt]
    \beta &\equiv \{i : y^*_i = 0 = w^*_i\} = \{1\}, \\[4pt]
    \gamma &\equiv \{i : y^*_i = 0 < w^*_i\} = \emptyset.
\end{align*}
Because the index set $\beta$ is not empty ($\beta = \{1\}$), the point lacks strict complementarity and is therefore degenerate. (Note: Although $F(1,0,0) = 1 \neq 0$, evaluating the local structure of the mapping depends strictly on these spatial coordinates.)

(b) Constructing the Piecewise Linear Mapping $LH^*$

First, we find the derivatives of the constraint $F(x,y,w) = x + y - w$ with respect to the lower-level variables:
$DF_y = \nabla_y F = 1$, and $DF_w = \nabla_w F = -1$.
Since $\alpha$ and $\gamma$ are empty and there is no free variable $z$ in this problem, the mapping $LH^*$ simplifies strictly to the $\beta$ components. By definition from Section 6.1.1:

$$LH^*(dy_\beta, dw_\beta) = \begin{pmatrix} DF_{y_\beta} dy_\beta + DF_{w_\beta} dw_\beta \\ \min(dy_\beta, dw_\beta) \end{pmatrix}$$
Substituting the derivatives:

$$LH^*(dy_\beta, dw_\beta) = \begin{pmatrix} dy_\beta - dw_\beta \\ \min(dy_\beta, dw_\beta) \end{pmatrix}$$
(c) The $\mathcal{W}$ Property and Global Homeomorphism

To use Lemma 6.1.1, we must construct the matrix pair $(W_0, W_1)$ associated with the linear pieces of $LH^*$:
$W_0 = [DF_{y_\beta}] = [1]$, and
$W_1 = [-DF_{w_\beta}] = [-(-1)] = [1]$.
The $\mathcal{W}$ property requires that all column representatives of the pair $(W_0, W_1)$ have determinants of the same sign.

Since both $W_0$ and $W_1$ are $1 \times 1$ matrices with a value (and determinant) of $+1$, all possible column configurations yield a determinant of $1$. Because the sign is strictly positive in all cases, the pair $(W_0, W_1)$ has the $\mathcal{W}$ property. Therefore, by Lemma 6.1.1, $LH^*$ is a global Lipschitzian homeomorphism.

\subsubsection*{Problem 2: PIPA Quadratic Subproblem Setup}
(a) Penalty Constraint Residual $\phi_0$

First, map the constraint $w = -2 + x + y$ to the standard form $r(x,y,w) = q + Nx + My - w$:
$$r(x,y,w) = -2 + x + y - w$$
At the initial iterate $x^0 = 1, y^0 = 2, w^0 = 1$, the residual is:
$$r^0 = -2 + 1 + 2 - 1 = 0$$
The penalty constraint residual $\phi_0$ is defined as:
$$\phi_0 = (y^0)^T w^0 + \|r^0\| = 2$$

(b) Explicit $QP_0$ Formulation

We need to calculate the gradient values and matrix coefficients at the current iterate.
The objective gradients are
\begin{align*}
    df_x^0 &= \nabla_x f(1,2) = x^0 + y^0 = 1 + 2 = 3, \\[4pt]
    df_y^0 &= \nabla_y f(1,2) = x^0 + 2y^0 = 1 + 4 = 5.
\end{align*}
The centrality/complementarity terms are:
\begin{align*}
    \mu_0 &= \frac{(y^0)^T w^0}{m} = \frac{(2)(1)}{1} = 2, \\[4pt]
    Y_0 &= \text{diag}(y^0) = [2], \\[4pt]
    W_0 &= \text{diag}(w^0) = [1].
\end{align*}
The linearization constraint is $Ndx + Mdy - dw = -r^0$.
$(1)dx + (1)dy - dw = 0$ implies $dx + dy - dw = 0$.

The centering constraint is $W_0 dy + Y_0 dw = -Y_0 w^0 + \sigma_0 \mu_0 e$.
$(1)dy + (2)dw = -(2)(1) + (0.5)(2) = -2 + 1 = -1$ implies $dy + 2dw = -1$

Putting it all together, the full $QP_0$ to solve for the search direction is:
\begin{align*}
    \text{minimize} \quad & 3dx + 5dy + \frac{1}{2}dx^2 \\
    \text{subject to} \quad & dx \ge -1 \\
    & dx^2 \le 20 \\
    & dx + dy - dw = 0 \\
    & dy + 2dw = -1
\end{align*}

\subsubsection*{Problem 3: Implicit Programming and Piecewise Subproblems} 

(a) Directional Derivative $d\overline{y}$

The lower level problem is equivalent to the LCP: $y \ge 0, \; y - x \ge 0, \; y(y-x) = 0$.
This means the implicit solution mapping is exactly the $\max$ function: $\overline{y}(x) = \max(0, x)$.
The directional derivative of $\overline{y}$ at $x^\nu = 0$ in the direction $dx$ is calculated by definition:

$$(\overline{y}^\nu)'(0; dx) = \lim_{\tau \to 0^+} \frac{\max(0, 0 + \tau dx) - \max(0, 0)}{\tau} = \lim_{\tau \to 0^+} \frac{\tau \max(0, dx)}{\tau} = \max(0, dx)$$.

(b) Formulating the Descent Subproblem

Evaluate the objective gradients at $(x^\nu, \overline{y}^\nu) = (0,0)$:
\begin{align*}
    df_x = \nabla_x (x^2 - y) = 2x = 0, \\[4pt]
    df_y = \nabla_y (x^2 - y) = -1.
\end{align*}
Substitute these into the direction-finding subproblem (Equation 6.3.2):
\begin{align*}
    \text{minimize} \quad & -\max(0, dx) + dx^2 \\
    \text{subject to}\quad & -1 \le dx \le 1
\end{align*}

(c) Solving the Subproblem and Stationarity Implication

We analyze the piecewise objective function $g(dx) = -\max(0, dx) + dx^2$ over the feasible interval $[-1, 1]$:
\begin{enumerate}[label=(\roman*)]
    \item Case 1 ($dx \in [-1, 0]$): $g(dx) = dx^2$. The minimum occurs at $dx = 0$, yielding an objective value of $0$.
    \item Case 2 ($dx \in [0, 1]$):
$g(dx) = -dx + dx^2$. We find the minimum by taking the derivative and setting it to zero:
$-1 + 2dx = 0 \implies dx = 0.5$.
The objective value at this point is $-(0.5) + (0.5)^2 = -0.5 + 0.25 = -0.25$.
\end{enumerate}
Comparing the cases, the globally optimal search direction is $dx^ = 0.5$*, which yields a strictly negative objective value of $-0.25$.

Implication: According to Proposition 6.3.2, an optimum objective value of zero is a necessary and sufficient condition for the current point to be a stationary point. Because our optimal solution yields a strictly negative objective value, $dx = 0$ is not optimal, proving that $x^\nu = 0$ is NOT a stationary point of the MPEC.

\subsection{Solutions for proof questions}

\subsubsection*{Problem 1: Matrix Properties in Complementarity}

\textbf{(a) Definition of the Mixed P Property} \\
According to Definition 6.1.4, a partitioned matrix $Q = \begin{bmatrix} A & B & C \end{bmatrix}$ where $A, B \in \mathbb{R}^{(m+l) \times m}$ and $C \in \mathbb{R}^{(m+l) \times l}$ is said to have the \textit{mixed P property} if $C$ has full column rank and the following implication holds:
\[
\left. \begin{array}{l}
A r + B s + C t = 0 \\
(r, s) \neq 0
\end{array} \right\} \implies r_i s_i > 0 \text{ for some } i \in \{1, \dots, m\}.
\]

\textbf{(b) Proof of Linear Independence (Proposition 6.1.5)} \\
Let $\alpha \subseteq \{1, \dots, m\}$ be an arbitrary index set and $\overline{\alpha}$ be its complement. We must show that the matrix $M = \begin{bmatrix} A_\alpha & B_{\overline{\alpha}} & C \end{bmatrix}$ has linearly independent columns. 

Assume, for the sake of contradiction, that the columns are linearly dependent. Then there exists a non-zero vector $(r_\alpha, s_{\overline{\alpha}}, t) \neq 0$ such that:
\[
A_\alpha r_\alpha + B_{\overline{\alpha}} s_{\overline{\alpha}} + C t = 0.
\]
We define full $m$-dimensional vectors $r$ and $s$ by padding with zeros:
\[
r \equiv (r_\alpha, 0_{\overline{\alpha}}) \quad \text{and} \quad s \equiv (0_\alpha, s_{\overline{\alpha}}).
\]
Substituting these into our system yields:
\[
A r + B s + C t = 0.
\]
Now, observe the Hadamard (element-wise) product $r \circ s$. For any index $i$, if $i \in \alpha$, then $s_i = 0$. If $i \in \overline{\alpha}$, then $r_i = 0$. In all cases, $r_i s_i = 0$. Therefore, $r \circ s = \mathbf{0}$.

Because $r_i s_i = 0$ for all $i$, there does \textit{not} exist any index $i$ such that $r_i s_i > 0$. By the contrapositive of the mixed P property definition, this strictly requires that $(r, s) = 0$. 

If $(r, s) = 0$, then $r_\alpha = 0$ and $s_{\overline{\alpha}} = 0$. Our system simplifies to:
\[
C t = 0.
\]
Since $Q$ having the mixed P property fundamentally requires $C$ to have full column rank, the null space of $C$ is trivial, meaning $t = 0$. 
Thus, $(r_\alpha, s_{\overline{\alpha}}, t) = 0$, which contradicts our initial assumption. Therefore, the columns of $\begin{bmatrix} A_\alpha & B_{\overline{\alpha}} & C \end{bmatrix}$ are linearly independent. $\hfill \blacksquare$

\subsubsection*{Problem 2: PIPA Merit Function and Descent}

\textbf{(a) Derivation of the Directional Derivative (Lemma 6.1.9)} \\
First, we compute the exact partial derivatives of the constraint violation function $\phi(x,y,w,z) = F(x,y,w,z)^T F(x,y,w,z) + w^T y$. Applying the chain rule, we obtain the gradients with respect to each variable block:
\begin{align*}
    d\phi_x &= \nabla_x \phi = 2 (\nabla_x F)^T F \\
    d\phi_y &= \nabla_y \phi = 2 (\nabla_y F)^T F + w \\
    d\phi_w &= \nabla_w \phi = 2 (\nabla_w F)^T F + y \\
    d\phi_z &= \nabla_z \phi = 2 (\nabla_z F)^T F
\end{align*}
Next, we substitute these into the directional derivative expression:
\begin{align*}
    \Delta \phi &= (d\phi_x)^T dx + (d\phi_y)^T dy + (d\phi_w)^T dw + (d\phi_z)^T dz \\
    &= 2 F^T (\nabla_x F dx + \nabla_y F dy + \nabla_w F dw + \nabla_z F dz) + w^T dy + y^T dw
\end{align*}
From the equality constraints of the QP subproblem, the Newton step for $F$ satisfies the exact linearization:
\[
\nabla_x F dx + \nabla_y F dy + \nabla_w F dw + \nabla_z F dz = -F
\]
Substituting this into our equation yields:
\[
\Delta \phi = 2 F^T (-F) + w^T dy + y^T dw = -2\|F\|^2 + \sum_{j=1}^m (w_j dy_j + y_j dw_j)
\]
From the perturbed complementarity constraints in the QP subproblem, we know:
\[
w_j dy_j + y_j dw_j = -w_j y_j + \sigma_\nu \frac{y^T w}{m}
\]
Summing this relation over all $j = 1, \dots, m$:
\begin{align*}
    \sum_{j=1}^m (w_j dy_j + y_j dw_j) &= -\sum_{j=1}^m w_j y_j + \sum_{j=1}^m \sigma_\nu \frac{y^T w}{m} \\
    &= -y^T w + m \left( \sigma_\nu \frac{y^T w}{m} \right) \\
    &= -y^T w + \sigma_\nu y^T w \\
    &= -(1 - \sigma_\nu) y^T w
\end{align*}
Substituting this final sum back into $\Delta \phi$:
\[
(d\phi_x^\nu)^T dx^\nu + (d\phi_y^\nu)^T dy^\nu + (d\phi_w^\nu)^T dw^\nu + (d\phi_z^\nu)^T dz^\nu = -2\|F^\nu\|^2 - (1 - \sigma_\nu)(y^\nu)^T w^\nu. \hfill \blacksquare
\]

\subsubsection*{Problem 3: Implicit Programming and Stationarity}

\textbf{(a) Optimum Objective Value of Zero} \\
Let the objective function of the subproblem be denoted as $\Phi(dx) = (df_x^\nu)^T dx + (df_y^\nu)^T (\overline{y}^\nu)'(x^\nu; dx) + \frac{1}{2} dx^T Q_\nu dx$.
Clearly, $dx = 0$ is a feasible solution (since $x^\nu \in X$), and $\Phi(0) = 0$. 

\textit{($\Rightarrow$)} Suppose the optimum objective value is zero. Since $Q_\nu$ is positive definite, $\frac{1}{2} dx^T Q_\nu dx > 0$ for all $dx \neq 0$. If there existed another global optimal solution $dx^* \neq 0$, it must yield $\Phi(dx^*) = 0$. This implies $(df_x^\nu)^T dx^* + (df_y^\nu)^T (\overline{y}^\nu)'(x^\nu; dx^*) = -\frac{1}{2} (dx^*)^T Q_\nu dx^* < 0$. 

By the positive homogeneity of the directional derivative (i.e., $g'(\tau d) = \tau g'(d)$ for $\tau > 0$), scaling $dx^*$ by some $\tau \in (0, 1)$ yields:
\[
\Phi(\tau dx^*) = \tau \left[ (df_x^\nu)^T dx^* + (df_y^\nu)^T (\overline{y}^\nu)'(x^\nu; dx^*) \right] + \frac{\tau^2}{2} (dx^*)^T Q_\nu dx^*
\]
Since the linear/directional term is strictly negative, for $\tau > 0$ sufficiently small, the $\mathcal{O}(\tau)$ term dominates the $\mathcal{O}(\tau^2)$ term, leading to $\Phi(\tau dx^*) < 0$. This contradicts the assumption that the global optimum value is $0$. Thus, $dx = 0$ is the unique globally optimal solution.

\textit{($\Leftarrow$)} Suppose $dx = 0$ is a globally optimal solution. Since $\Phi(0) = 0$, the optimum objective value of the subproblem is $0$. $\hfill \blacksquare$

\vspace{0.5cm}
\textbf{(b) Stationarity Proof} \\
Assume the optimum objective value of the subproblem is $0$. Then, for any feasible direction $dx \in X - x^\nu$ and any scalar $\tau \in (0, 1]$, the vector $x^\nu + \tau dx \in X$ (by the convexity of the polyhedron $X$). Since the minimum is $0$, we evaluate the objective at the feasible step $\tau dx$:
\[
\tau(df_x^\nu)^T dx + \tau(df_y^\nu)^T (\overline{y}^\nu)'(x^\nu; dx) + \frac{\tau^2}{2} dx^T Q_\nu dx \ge 0.
\]
Dividing the entire inequality by $\tau > 0$ yields:
\[
(df_x^\nu)^T dx + (df_y^\nu)^T (\overline{y}^\nu)'(x^\nu; dx) + \frac{\tau}{2} dx^T Q_\nu dx \ge 0.
\]
Taking the limit as $\tau \to 0^+$, the quadratic term vanishes:
\[
(df_x^\nu)^T dx + (df_y^\nu)^T (\overline{y}^\nu)'(x^\nu; dx) \ge 0.
\]
Because $X$ is a polyhedron, the tangent cone $\mathcal{T}(x^\nu; X)$ is equivalent to the closure of the cone of feasible directions, $\bigcup_{\tau > 0} \tau(X - x^\nu)$. Therefore, the inequality holds for all $dx \in \mathcal{T}(x^\nu; X)$. By Lemma 4.2.5, this variational inequality over the tangent cone is the precise definition of stationarity for the implicit programming formulation at $(x^\nu, y^\nu)$. $\hfill \blacksquare$

\subsubsection*{Problem 4: Polyhedral Geometry in Piecewise SQP}

\textbf{(a) Proof of Lemma 6.3.3} \\
Represent the polyhedron $X$ using a finite set of linear inequalities: $X = \{x \in \mathbb{R}^n : A x \le b\}$. Let $\mathcal{I}(\overline{x}) = \{i : A_i \overline{x} = b_i\}$ be the index set of active (binding) constraints at $\overline{x}$. 

For a polyhedron, the tangent cone at $\overline{x}$ is exactly characterized by the active constraints:
\[
\mathcal{T}(\overline{x}; X) = \{d \in \mathbb{R}^n : A_i d \le 0, \;\; \forall i \in \mathcal{I}(\overline{x})\}.
\]
Let $d \in \mathcal{T}(\overline{x}; X)$. 

1. \textbf{Active Constraints ($i \in \mathcal{I}(\overline{x})$):} 
For any $y \in X$ (which means $A_i y \le b_i$) and any $\epsilon \ge 0$:
\[
A_i(y + \epsilon d) = A_i y + \epsilon A_i d.
\]
Since $y \in X$, $A_i y \le b_i$, and since $d \in \mathcal{T}(\overline{x}; X)$, $A_i d \le 0$. Thus, $A_i(y + \epsilon d) \le b_i + 0 = b_i$. This holds for all $\epsilon \ge 0$.

2. \textbf{Inactive Constraints ($i \notin \mathcal{I}(\overline{x})$):}
For these constraints, $A_i \overline{x} < b_i$. By continuity, there exists a $\delta > 0$ and $c > 0$ such that for any $y$ satisfying $\|y - \overline{x}\| \le \delta$, we have $A_i y \le b_i - c$. 
We must find an $\overline{\epsilon} > 0$ small enough such that the step $\epsilon A_i d$ does not exceed the slack $c$. Define:
\[
\overline{\epsilon} = \min_{i \notin \mathcal{I}(\overline{x}), A_i d > 0} \frac{c}{A_i d}
\]
(If $A_i d \le 0$ for all inactive constraints, any $\overline{\epsilon} > 0$ works). 
Now, for any $\epsilon \in [0, \overline{\epsilon}]$ and any $y \in X$ such that $\|y - \overline{x}\| \le \delta$:
\[
A_i(y + \epsilon d) = A_i y + \epsilon A_i d \le (b_i - c) + c = b_i.
\]
Since $A_i(y + \epsilon d) \le b_i$ holds for both active and inactive constraints under these conditions, $y + \epsilon d \in X$. $\hfill \blacksquare$

\textbf{(b) Application in Limit Transitions (Theorem 6.3.4)} \\
In the convergence proof of the descent algorithm, we must show that the limit point $x^*$ satisfies the stationarity variational inequality over its tangent cone: $(df_x^*)^T dx + (df_y^*)^T (\overline{y}^*)'(x^*; dx) \ge 0$ for all $dx \in \mathcal{T}(x^*; X)$. The fundamental difficulty is that while $dx$ is an admissible direction at the limit $x^*$, it may point slightly \textit{outside} the feasible set from the perspective of an intermediate iterate $x^\nu \neq x^*$. 

Lemma 6.3.3 resolves this by ensuring a uniform geometric cushion. By setting $y = x^\nu$ and $\overline{x} = x^*$, the lemma guarantees that once the sequence converges close enough to the limit ($\|x^\nu - x^*\| \le \delta$), a finite, uniform step $\overline{\epsilon} > 0$ exists such that $x^\nu + \epsilon dx$ remains strictly inside $X$ for all subsequent $\nu$. This permits the construction of a valid feasible test step $\epsilon dx$ across the entire tail of the sequence, allowing the $\liminf$ operator to securely pass through the descent bounds to establish stationarity at the limit without feasibility violations blocking the path.

\newpage
\section{Takeaways}

\begin{enumerate}[label=(\arabic*),itemsep=0.3em]
    \item MPEC algorithms are difficult because complementarity produces a feasible set that is
    piecewise smooth rather than globally smooth.
    \item PIPA is best understood as a synthesis of SQP, Newton linearization, and
    central-path regularization.
    \item Leyffer's counterexample is not merely a technicality: it identifies a genuine
    design issue in the coupling between upper-level motion and lower-level feasibility.
    \item The monotone-LCP variant repairs this coupling by preventing complementarity from
    vanishing too quickly relative to the affine residual.
    \item Implicit-programming methods are reduced-space methods; they are mathematically
    elegant but often limited by the cost of evaluating the lower-level response map.
    \item PSQP shows that local fast convergence can be justified directly on the piecewise
    geometry of the MPEC feasible set.
\end{enumerate}

\section{Bibliographic remarks and Acknowledgment}

For the classical MPEC framework and the original PIPA development, see the monograph
of Luo, Pang, and Ralph. For the failure of the strongest original PIPA convergence claim,
see Leyffer's counterexample. For the piecewise SQP viewpoint, see Ralph's PSQP treatment. 
Please see the relevant references: \cite{
leyffer_penalty_2005,
facchinei_ghost_2021,
kuhn_piecewise_1987,
cottle_linear_1992,
kojima_unified_1991,
liu_bidirectional_2025,
larson_manifold_2021,
hu_convergence_2026,
todd_centered_1990,
wright_superlinear_1996,
zhang_convergence_1994,
jm_iterative_1970,
burke_sequential_1989,
aghasi_fully_2025,
burke_robust_1992,
wang_analysis_2025,
lampariello_detecting_2025,
londe_multi-stage_2025,
hu_improved_2023,
boob_level_2025,
yao_relative_2023,
dussault_polyhedral_2026,
caselli_bilevel_2026,
burke_robust_1989,
wang_interior_1996,
el_hajj_efficient_2026,
pang_minimization_1991,
lu_first-order_2024,
chen_two_2025,
lewis_activeset_2021,
wang_algebraicspectral_2026,
lagos_bilevel_2025,
cui_complexity_2023,
casas_convergence_2024,
outrata_numerical_1995,
mccormick_nonlinear_1983,
bonnans_local_1994,
pang_convergence_1993,
liang_global_2025,
malitsky_projected_2015,
kyparisis_uniqueness_1985,
robinson_perturbed_1974,
wang2026lecture,
curtis_eaves_where_1983,
wang_analysis_2025-1,
han_continuous_2024,
pang_nesqp_1993,
bonnans_local_1989,
gao2022rolling,
wang2025multi,
christof_nonsmooth_2020,
han_globally_1977,
lou_design_2025,
burke_inexact_2020,
wang_damage-structured_2026,
beck_survey_2023,
dai_fixed_2020,
gowda_stability_1994,
shin_exponential_2022,
obara_sequential_2022,
yu_pattern_2026,
jiang_extension_2003,
dempe_necessary_1992,
hansen_new_1992,
wang_sequential_2023,
rawat_augmented_2026,
mordukhovich_generalized_2021,
mohammadi_variational_2022,
duy_khanh_generalized_2023,
izmailov_perturbed_2022}. 
\textbf{This note is mainly based on Chapter 6,
\emph{First-Order Optimality Conditions}, in the MPEC monograph of Zhi-Quan Luo, Jong-Shi Pang and Daniel Ralph.}

\bibliographystyle{unsrtnat}
\bibliography{references}

@article{kuhn_piecewise_1987,
	title = {Piecewise affine bijections of {Rn}, and the equation {Sx}+-{Tx}-=y},
	volume = {96},
	copyright = {https://www.elsevier.com/tdm/userlicense/1.0/},
	doi = {10.1016/0024-3795(87)90339-9},
	language = {en},
	journal = {Linear Algebra and its Applications},
	author = {Kuhn, D. and Löwen, R.},
	year = {1987},
	pages = {109--129},
}

@book{cottle_linear_1992,
	address = {Boston},
	edition = {1st},
	title = {The {Linear} {Complementarity} {Problem}},
	publisher = {Academic Press},
	author = {Cottle, R.W. and Pang, J.S. and Stone, R.E.},
	year = {1992},
}

@article{kojima_unified_1991,
	title = {A unified approach to interior point algorithms for linear complementarity problems: {A} summary},
	volume = {10},
	shorttitle = {A unified approach to interior point algorithms for linear complementarity problems},
	doi = {10.1016/0167-6377(91)90010-M},
	language = {en},
	number = {5},
	journal = {Operations Research Letters},
	author = {Kojima, M. and Megiddo, N. and Noma, T. and Yoshise, A.},
	year = {1991},
	pages = {247--254},
}

@article{todd_centered_1990,
	title = {A centered projective algorithm for linear programming},
	volume = {15},
	url = {https://www.jstor.org/stable/3689994},
	number = {3},
	urldate = {2026-04-06},
	journal = {Mathematics of Operations Research},
	publisher = {INFORMS},
	author = {Todd, M. J. and Ye, Y.},
	year = {1990},
	pages = {508--529},
}

@article{wright_superlinear_1996,
	title = {A superlinear infeasible-interior-point algorithm for monotone complementarity problems},
	volume = {21},
	url = {https://www.jstor.org/stable/3690189},
	number = {4},
	urldate = {2026-04-06},
	journal = {Mathematics of Operations Research},
	publisher = {INFORMS},
	author = {Wright, S. and Ralph, D.},
	year = {1996},
	pages = {815--838},
}

@article{zhang_convergence_1994,
	title = {On the convergence of a class of infeasible interior-point methods for the horizontal linear complementarity problem},
	volume = {4},
	doi = {10.1137/0804012},
	language = {en},
	number = {1},
	journal = {SIAM Journal on Optimization},
	author = {Zhang, Y.},
	year = {1994},
	pages = {208--227},
}

@book{jm_iterative_1970,
	address = {New York},
	title = {Iterative {Solution} of {Nonlinear} {Equations} in {Several} {Variables}},
	url = {https://www.scirp.org/reference/referencespapers?referenceid=1205027},
	urldate = {2026-04-06},
	publisher = {Academic Press},
	author = {J.M., , Ortega and W.G., , Rheinboldt},
	year = {1970},
}

@article{burke_sequential_1989,
	title = {A sequential quadratic programming method for potentially infeasible mathematical programs},
	volume = {139},
	copyright = {https://www.elsevier.com/tdm/userlicense/1.0/},
	doi = {10.1016/0022-247X(89)90111-X},
	language = {en},
	number = {2},
	journal = {Journal of Mathematical Analysis and Applications},
	author = {Burke, J. V.},
	year = {1989},
	pages = {319--351},
}

@article{burke_robust_1992,
	title = {A robust trust region method for constrained nonlinear programming problems},
	volume = {2},
	doi = {10.1137/0802016},
	language = {en},
	number = {2},
	journal = {SIAM Journal on Optimization},
	author = {Burke, J. V.},
	year = {1992},
	pages = {325--347},
}

@article{burke_robust_1989,
	title = {A robust sequential quadratic programming method},
	volume = {43},
	copyright = {http://www.springer.com/tdm},
	doi = {10.1007/BF01582294},
	language = {en},
	number = {1-3},
	journal = {Mathematical Programming},
	author = {Burke, J. V. and Han, S-P.},
	year = {1989},
	pages = {277--303},
}

@article{wang_interior_1996,
	title = {An interior point potential reduction method for constrained equations},
	volume = {74},
	copyright = {http://www.springer.com/tdm},
	doi = {10.1007/BF02592210},
	language = {en},
	number = {2},
	journal = {Mathematical Programming},
	author = {Wang, T. and Monteiro, R. D. C. and Pang, J.-S.},
	year = {1996},
	pages = {159--195},
}

@article{pang_minimization_1991,
	title = {Minimization of locally lipschitzian functions},
	volume = {1},
	doi = {10.1137/0801006},
	language = {en},
	number = {1},
	journal = {SIAM Journal on Optimization},
	author = {Pang, J.-S. and Han, S.-P. and Rangaraj, N.},
	year = {1991},
	pages = {57--82},
}

@article{outrata_numerical_1995,
	title = {A numerical approach to optimization problems with variational inequality constraints},
	volume = {68},
	copyright = {http://www.springer.com/tdm},
	doi = {10.1007/BF01585759},
	language = {en},
	number = {1-3},
	journal = {Mathematical Programming},
	author = {Outrata, J. and Zowe, J.},
	year = {1995},
	pages = {105--130},
}

@book{mccormick_nonlinear_1983,
	address = {New York},
	title = {Nonlinear programming: {Theory}, algorithms, and applications},
	isbn = {978-0-471-09309-1},
	shorttitle = {Nonlinear programming},
	publisher = {Wiley},
	author = {McCormick, G. P.},
	year = {1983},
}

@article{bonnans_local_1994,
	title = {Local analysis of {Newton}-type methods for variational inequalities and nonlinear programming},
	volume = {29},
	copyright = {http://www.springer.com/tdm},
	doi = {10.1007/BF01204181},
	language = {en},
	number = {2},
	journal = {Applied Mathematics \& Optimization},
	author = {Bonnans, J. F.},
	year = {1994},
	pages = {161--186},
}

@article{pang_convergence_1993,
	title = {Convergence of splitting and {Newton} methods for complementarity problems: {An} application of some sensitivity results},
	volume = {58},
	copyright = {http://www.springer.com/tdm},
	shorttitle = {Convergence of splitting and {Newton} methods for complementarity problems},
	doi = {10.1007/BF01581264},
	language = {en},
	number = {1-3},
	journal = {Mathematical Programming},
	author = {Pang, J.-S.},
	year = {1993},
	pages = {149--160},
}

@article{kyparisis_uniqueness_1985,
	title = {On uniqueness of {Kuhn}-{Tucker} multipliers in nonlinear programming},
	volume = {32},
	copyright = {http://www.springer.com/tdm},
	doi = {10.1007/BF01586095},
	language = {en},
	number = {2},
	journal = {Mathematical Programming},
	author = {Kyparisis, J.},
	year = {1985},
	pages = {242--246},
}

@article{robinson_perturbed_1974,
	title = {Perturbed {Kuhn}-{Tucker} points and rates of convergence for a class of nonlinear-programming algorithms},
	volume = {7},
	copyright = {http://www.springer.com/tdm},
	doi = {10.1007/BF01585500},
	language = {en},
	number = {1},
	journal = {Mathematical Programming},
	author = {Robinson, S. M.},
	year = {1974},
	pages = {1--16},
}

@incollection{curtis_eaves_where_1983,
	address = {Boston, MA},
	title = {Where {Solving} for {Stationary} {Points} by {LCPs} is {Mixing} {Newton} {Iterates}},
	isbn = {978-1-4613-3574-0 978-1-4613-3572-6},
	doi = {10.1007/978-1-4613-3572-6_5},
	language = {en},
	booktitle = {Homotopy {Methods} and {Global} {Convergence}},
	publisher = {Springer US},
	author = {Curtis Eaves, B.},
	year = {1983},
	pages = {63--77},
}

@article{pang_nesqp_1993,
	title = {{NE}/{SQP}: {A} robust algorithm for the nonlinear complementarity problem},
	volume = {60},
	copyright = {http://www.springer.com/tdm},
	shorttitle = {{NE}/{SQP}},
	doi = {10.1007/BF01580617},
	language = {en},
	number = {1-3},
	journal = {Mathematical Programming},
	author = {Pang, J.-S. and Gabriel, S. A.},
	year = {1993},
	pages = {295--337},
}

@incollection{bonnans_local_1989,
	title = {Local study of newton type algorithms for constrained problems},
	volume = {1405},
	isbn = {978-3-540-51970-6},
	doi = {10.1007/BFb0083583},
	language = {en},
	booktitle = {Optimization},
	publisher = {Springer Berlin Heidelberg},
	author = {Bonnans, J. F.},
	year = {1989},
	note = {Series Title: Lecture Notes in Mathematics},
	pages = {13--24},
}

@article{han_globally_1977,
	title = {A globally convergent method for nonlinear programming},
	volume = {22},
	copyright = {http://www.springer.com/tdm},
	doi = {10.1007/BF00932858},
	language = {en},
	number = {3},
	journal = {Journal of Optimization Theory and Applications},
	author = {Han, S. P.},
	year = {1977},
	pages = {297--309},
}

@article{gowda_stability_1994,
	title = {Stability analysis of variational inequalities and nonlinear complementarity problems, via the mixed linear complementarity problem and degree theory},
	volume = {19},
	doi = {10.1287/moor.19.4.831},
	language = {en},
	number = {4},
	journal = {Mathematics of Operations Research},
	author = {Gowda, M. S. and Pang, J.-S.},
	year = {1994},
	pages = {831--879},
}

@article{jiang_extension_2003,
	title = {Extension of quasi-newton methods to mathematical programs with complementarity constraints},
	volume = {25},
	copyright = {https://www.springernature.com/gp/researchers/text-and-data-mining},
	doi = {10.1023/A:1022945316191},
	language = {en},
	number = {1-3},
	journal = {Computational Optimization and Applications},
	author = {Jiang, H. and Ralph, D.},
	year = {2003},
	pages = {123--150},
}

@article{dempe_necessary_1992,
	title = {A necessary and a sufficient optimality condition for bilevel programming problems},
	volume = {25},
	doi = {10.1080/02331939208843831},
	language = {en},
	number = {4},
	journal = {Optimization},
	author = {Dempe, S.},
	year = {1992},
	pages = {341--354},
}

@article{hansen_new_1992,
	title = {New branch-and-bound rules for linear bilevel programming},
	volume = {13},
	doi = {10.1137/0913069},
	language = {en},
	number = {5},
	journal = {SIAM Journal on Scientific and Statistical Computing},
	author = {Hansen, P. and Jaumard, B. and Savard, G.},
	year = {1992},
	pages = {1194--1217},
}

@article{liu_bidirectional_2025,
	title = {Bidirectional endothelial feedback drives turing-vascular patterning and drug-resistance niches: a hybrid {PDE}-agent-based study},
	volume = {12},
	shorttitle = {Bidirectional endothelial feedback drives turing-vascular patterning and drug-resistance niches},
	doi = {10.3390/bioengineering12101097},
	language = {en},
	number = {10},
	journal = {Bioengineering},
	author = {Liu, Z. and Wang, L. S. and Yu, J. and Zhang, J. and Martel, E. and Li, S.},
	year = {2025},
	pages = {1097},
}

@article{wang_analysis_2025,
	title = {Analysis framework for stochastic predator–prey model with demographic noise},
	volume = {27},
	doi = {10.1016/j.rinam.2025.100621},
	language = {en},
	journal = {Results in Applied Mathematics},
	author = {Wang, L. S. and Yu, J.},
	year = {2025},
	pages = {100621},
}

@article{wang_algebraicspectral_2026,
	title = {Algebraic–spectral thresholds and discrete–continuous stability transfer in {Leslie}–{Gower} systems},
	volume = {34},
	doi = {10.3934/era.2026013},
	number = {1},
	journal = {Electronic Research Archive},
	author = {Wang, L. S. and Yu, J.},
	year = {2026},
	pages = {251--290},
}

@article{liang_global_2025,
	title = {Global well-posedness and stability of nonlocal damage-structured lineage model with feedback and dedifferentiation},
	volume = {13},
	doi = {10.3390/math13223583},
	language = {en},
	number = {22},
	journal = {Mathematics},
	author = {Liang, Y. and Wang, L. S. and Yu, J. and Liu, Z.},
	year = {2025},
	pages = {3583},
}

@article{wang_analysis_2025-1,
	title = {Analysis and mean-field limit of a hybrid {PDE}-{ABM} modeling angiogenesis-regulated resistance evolution},
	volume = {13},
	doi = {10.3390/math13172898},
	language = {en},
	number = {17},
	journal = {Mathematics},
	author = {Wang, L. S. and Yu, J. and Li, S. and Liu, Z.},
	year = {2025},
	pages = {2898},
}

@article{wang_damage-structured_2026,
	title = {A damage-structured {PDE} model of stem cell hierarchies: {The} dual role of dedifferentiation in tissue homeostasis and aging},
	volume = {21},
	shorttitle = {A damage-structured {PDE} model of stem cell hierarchies},
	doi = {10.1371/journal.pone.0335163},
	language = {en},
	number = {2},
	journal = {PLOS One},
	author = {Wang, L. S. and Yu, J. and Liu, Z.},
	year = {2026},
	pages = {e0335163},
}

@article{yu_pattern_2026,
	title = {Pattern suppression and recovery under one-way versus two-way chemotactic coupling in hybrid partial differential equation–ordinary differential equation models},
	issn = {3052-878X},
	doi = {10.1515/tp-2026-0023},
	language = {en},
	journal = {Transport Phenomena},
	author = {Yu, J. and Wang, L. S. and Liu, Z. and Liu, J.},
	year = {2026},
}

@article{casas_convergence_2024,
	title = {Convergence analysis of the semismooth newton method for sparse control problems governed by semilinear elliptic equations},
	volume = {62},
	doi = {10.1137/23M1585945},
	language = {en},
	number = {6},
	journal = {SIAM Journal on Control and Optimization},
	author = {Casas, E. and Mateos, M.},
	year = {2024},
	pages = {3076--3090},
}

@article{lou_design_2025,
	title = {Design guidelines for noise-tolerant optimization with applications in robust design},
	volume = {47},
	doi = {10.1137/24M1632279},
	language = {en},
	number = {3},
	journal = {SIAM Journal on Scientific Computing},
	author = {Lou, Y. and Sun, S. and Nocedal, J.},
	year = {2025},
	pages = {A1335--A1357},
}

@article{hu_convergence_2026,
	title = {Convergence properties of gradient-based methods for minimax problems with nonlinear constraints},
	volume = {209},
	doi = {10.1007/s10957-026-02964-w},
	language = {en},
	number = {1},
	journal = {Journal of Optimization Theory and Applications},
	author = {Hu, Q. and Wang, B. and Xu, M.},
	year = {2026},
	pages = {13},
}

@article{shin_exponential_2022,
	title = {Exponential decay of sensitivity in graph-structured nonlinear programs},
	volume = {32},
	doi = {10.1137/21M1391079},
	language = {en},
	number = {2},
	journal = {SIAM Journal on Optimization},
	author = {Shin, S. and Anitescu, M. and Zavala, V. M.},
	year = {2022},
	pages = {1156--1183},
}

@article{malitsky_projected_2015,
	title = {Projected reflected gradient methods for monotone variational inequalities},
	volume = {25},
	doi = {10.1137/14097238X},
	language = {en},
	number = {1},
	journal = {SIAM Journal on Optimization},
	author = {Malitsky, Y.},
	year = {2015},
	pages = {502--520},
}

@article{dussault_polyhedral_2026,
	title = {Polyhedral newton-min algorithms for complementarity problems},
	volume = {215},
	doi = {10.1007/s10107-025-02219-y},
	language = {en},
	number = {1-2},
	journal = {Mathematical Programming},
	author = {Dussault, J.-P. and Frappier, M. and Gilbert, J. C.},
	year = {2026},
	pages = {269--324},
}

@article{dai_fixed_2020,
	title = {A fixed point iterative method for tensor complementarity problems},
	volume = {84},
	doi = {10.1007/s10915-020-01299-6},
	language = {en},
	number = {3},
	journal = {Journal of Scientific Computing},
	author = {Dai, P.-F.},
	year = {2020},
	pages = {49},
}

@article{yao_relative_2023,
	title = {Relative {Lipschitz}-like property of parametric systems via projectional coderivatives},
	volume = {33},
	doi = {10.1137/22M151296X},
	language = {en},
	number = {3},
	journal = {SIAM Journal on Optimization},
	author = {Yao, W. and Yang, X.},
	year = {2023},
	pages = {2021--2040},
}

@article{han_continuous_2024,
	title = {Continuous selections of solutions to parametric variational inequalities},
	volume = {34},
	doi = {10.1137/22M1514982},
	language = {en},
	number = {1},
	journal = {SIAM Journal on Optimization},
	author = {Han, S. and Pang, J.-S.},
	year = {2024},
	pages = {870--892},
}

@article{beck_survey_2023,
	title = {A survey on bilevel optimization under uncertainty},
	volume = {311},
	doi = {10.1016/j.ejor.2023.01.008},
	language = {en},
	number = {2},
	journal = {European Journal of Operational Research},
	author = {Beck, Y. and Ljubic, I. and Schmidt, M.},
	year = {2023},
	pages = {401--426},
}

@article{caselli_bilevel_2026,
	title = {Bilevel optimization with sustainability perspective: {A} survey on applications},
	volume = {332},
	shorttitle = {Bilevel optimization with sustainability perspective},
	doi = {10.1016/j.ejor.2025.08.051},
	language = {en},
	number = {2},
	journal = {European Journal of Operational Research},
	author = {Caselli, G. and Iori, M. and Ljubic, I.},
	year = {2026},
	pages = {357--375},
}

@article{lu_first-order_2024,
	title = {First-order penalty methods for bilevel optimization},
	volume = {34},
	doi = {10.1137/23M1566753},
	language = {en},
	number = {2},
	journal = {SIAM Journal on Optimization},
	author = {Lu, Z. and Mei, S.},
	year = {2024},
	pages = {1937--1969},
}

@article{lagos_bilevel_2025,
	title = {Bilevel optimization approach for fuel treatment planning},
	volume = {320},
	doi = {10.1016/j.ejor.2024.07.014},
	language = {en},
	number = {1},
	journal = {European Journal of Operational Research},
	author = {Lagos, T. and Choi, J. and Segundo, B. and Gan, J. and Ntaimo, L. and Prokopyev, O. A.},
	year = {2025},
	pages = {205--218},
}

@article{aghasi_fully_2025,
	title = {Fully zeroth-order bilevel programming via gaussian smoothing},
	volume = {205},
	doi = {10.1007/s10957-025-02647-y},
	language = {en},
	number = {2},
	journal = {Journal of Optimization Theory and Applications},
	author = {Aghasi, A. and Ghadimi, S.},
	year = {2025},
	pages = {31},
}

@article{londe_multi-stage_2025,
	title = {A multi-stage approach for {Root} {Sequence} {Index} allocation},
	volume = {327},
	doi = {10.1016/j.ejor.2025.05.015},
	language = {en},
	number = {1},
	journal = {European Journal of Operational Research},
	author = {Londe, M. A. and Andrade, C. E. and Pessoa, L. S.},
	year = {2025},
	pages = {95--114},
}

@article{hu_improved_2023,
	title = {An improved unconstrained approach for bilevel optimization},
	volume = {33},
	doi = {10.1137/22M1513034},
	language = {en},
	number = {4},
	journal = {SIAM Journal on Optimization},
	author = {Hu, X. and Xiao, N. and Liu, X. and Toh, K.-C.},
	year = {2023},
	pages = {2801--2829},
}

@article{el_hajj_efficient_2026,
	title = {Efficient resource sharing for strategic disaster preparedness},
	volume = {329},
	doi = {10.1016/j.ejor.2025.08.025},
	language = {en},
	number = {3},
	journal = {European Journal of Operational Research},
	author = {El Hajj, H. and Elhedhli, S. and Gzara, F.},
	year = {2026},
	pages = {836--847},
}

@article{boob_level_2025,
	title = {Level constrained first order methods for function constrained optimization},
	volume = {209},
	doi = {10.1007/s10107-024-02057-4},
	language = {en},
	number = {1-2},
	journal = {Mathematical Programming},
	author = {Boob, D. and Deng, Q. and Lan, G.},
	year = {2025},
	pages = {1--61},
}

@article{wang_sequential_2023,
	title = {A sequential quadratic programming algorithm for nonsmooth problems with upper-{C2} objective},
	volume = {33},
	doi = {10.1137/22M1490995},
	language = {en},
	number = {3},
	journal = {SIAM Journal on Optimization},
	author = {Wang, J. and Petra, C. G.},
	year = {2023},
	pages = {2379--2405},
}

@article{facchinei_ghost_2021,
	title = {Ghost penalties in nonconvex constrained optimization: {Diminishing} stepsizes and iteration complexity},
	volume = {46},
	shorttitle = {Ghost penalties in nonconvex constrained optimization},
	doi = {10.1287/moor.2020.1079},
	language = {en},
	number = {2},
	journal = {Mathematics of Operations Research},
	author = {Facchinei, F. and Kungurtsev, V. and Lampariello, L. and Scutari, G.},
	year = {2021},
	pages = {595--627},
}

@article{obara_sequential_2022,
	title = {Sequential quadratic optimization for nonlinear optimization problems on riemannian manifolds},
	volume = {32},
	doi = {10.1137/20M1370173},
	language = {en},
	number = {2},
	journal = {SIAM Journal on Optimization},
	author = {Obara, M. and Okuno, T. and Takeda, A.},
	year = {2022},
	pages = {822--853},
}

@article{burke_inexact_2020,
	title = {Inexact sequential quadratic optimization with penalty parameter updates within the {QP} solver},
	volume = {30},
	doi = {10.1137/18M1176488},
	language = {en},
	number = {3},
	journal = {SIAM Journal on Optimization},
	author = {Burke, J. V. and Curtis, F. E. and Wang, H. and Wang, J.},
	year = {2020},
	pages = {1822--1849},
}

@article{lampariello_detecting_2025,
	title = {On detecting degenerate stationarity},
	volume = {206},
	issn = {0022-3239, 1573-2878},
	doi = {10.1007/s10957-025-02747-9},
	language = {en},
	number = {3},
	journal = {Journal of Optimization Theory and Applications},
	author = {Lampariello, L.},
	year = {2025},
	pages = {69},
}

@article{larson_manifold_2021,
	title = {Manifold sampling for optimizing nonsmooth nonconvex compositions},
	volume = {31},
	doi = {10.1137/20M1378089},
	language = {en},
	number = {4},
	journal = {SIAM Journal on Optimization},
	author = {Larson, J. and Menickelly, M. and Zhou, B.},
	year = {2021},
	pages = {2638--2664},
}

@article{cui_complexity_2023,
	title = {Complexity guarantees for an implicit smoothing-enabled method for stochastic {MPECs}},
	volume = {198},
	doi = {10.1007/s10107-022-01893-6},
	language = {en},
	number = {2},
	journal = {Mathematical Programming},
	author = {Cui, S. and Shanbhag, U. V. and Yousefian, F.},
	year = {2023},
	pages = {1153--1225},
}

@article{christof_nonsmooth_2020,
	title = {A nonsmooth trust-region method for locally lipschitz functions with application to optimization problems constrained by variational inequalities},
	volume = {30},
	doi = {10.1137/18M1164925},
	language = {en},
	number = {3},
	journal = {SIAM Journal on Optimization},
	author = {Christof, C. and De Los Reyes, J. C. and Meyer, C.},
	year = {2020},
	pages = {2163--2196},
}

@article{rawat_augmented_2026,
	title = {Augmented {Lagrangian} neural network for solving mathematical programs with equilibrium constraints},
	volume = {209},
	doi = {10.1007/s10957-026-02963-x},
	language = {en},
	number = {1},
	journal = {Journal of Optimization Theory and Applications},
	author = {Rawat, A. and Singh, V.},
	year = {2026},
	pages = {15},
}

@article{mordukhovich_generalized_2021,
	title = {Generalized newton algorithms for tilt-stable minimizers in nonsmooth optimization},
	volume = {31},
	doi = {10.1137/20M1329937},
	language = {en},
	number = {2},
	journal = {SIAM Journal on Optimization},
	author = {Mordukhovich, B. S. and Sarabi, M. E.},
	year = {2021},
	pages = {1184--1214},
}

@article{mohammadi_variational_2022,
	title = {Variational analysis of composite models with applications to continuous optimization},
	volume = {47},
	doi = {10.1287/moor.2020.1074},
	language = {en},
	number = {1},
	journal = {Mathematics of Operations Research},
	author = {Mohammadi, A. and Mordukhovich, B. S. and Sarabi, M. E.},
	year = {2022},
	pages = {397--426},
}

@article{duy_khanh_generalized_2023,
	title = {A generalized newton method for subgradient systems},
	volume = {48},
	doi = {10.1287/moor.2022.1320},
	language = {en},
	number = {4},
	journal = {Mathematics of Operations Research},
	author = {Duy Khanh, P. and Mordukhovich, B. and Phat, V. T.},
	year = {2023},
	pages = {1811--1845},
}

@article{chen_two_2025,
	title = {Two typical implementable semismooth* {Newton} methods for generalized equations are {G}-semismooth newton methods},
	doi = {10.1287/moor.2024.0617},
	language = {en},
	journal = {Mathematics of Operations Research},
	author = {Chen, L. and Sun, D. and Zhang, W.},
	year = {2025},
	pages = {moor.2024.0617},
}

@article{lewis_activeset_2021,
	title = {Active‐set {Newton} methods and partial smoothness},
	volume = {46},
	doi = {10.1287/moor.2020.1075},
	language = {en},
	number = {2},
	journal = {Mathematics of Operations Research},
	author = {Lewis, A. S. and Wylie, C.},
	year = {2021},
	pages = {712--725},
}

@article{izmailov_perturbed_2022,
	title = {Perturbed augmented lagrangian method framework with applications to proximal and smoothed variants},
	volume = {193},
	doi = {10.1007/s10957-021-01914-y},
	language = {en},
	number = {1-3},
	journal = {Journal of Optimization Theory and Applications},
	author = {Izmailov, A. F. and Solodov, M. V.},
	year = {2022},
	pages = {491--522},
}

@article{leyffer_penalty_2005,
	title = {The penalty interior-point method fails to converge},
	volume = {20},
	doi = {10.1080/10556780500140078},
	language = {en},
	number = {4-5},
	journal = {Optimization Methods and Software},
	author = {Leyffer, S.},
	year = {2005},
	pages = {559--568},
}

@inproceedings{gao2022rolling,
  title={Rolling prediction model of closing price based on EEMD data noise reduction and HGS-DELM},
  author={Gao, Yuansheng and Li, Lei and Yu, Jiguang},
  booktitle={2022 International Conference on Data Analytics, Computing and Artificial Intelligence (ICDACAI)},
  pages={255--260},
  year={2022},
  organization={IEEE}
}

@article{wang2025multi,
  title={Multi-strategy Hybrid Improved Intelligent Algorithm for Solving UAV-MTSP},
  author={Wang, Zixin and Wang, Danqing and Yu, Jiguang},
  journal={Information Technology and Control},
  volume={54},
  number={2},
  pages={413--438},
  year={2025}
}

@article{wang2026lecture,
  title={Lecture Note for Bounded Controls in Continuous-Time and Control of Several Variables},
  author={Wang, L. S.},
  journal={arXiv preprint arXiv:2604.05882},
  year={2026}
}

\end{document}